\documentclass[11 pt]{amsart}

\usepackage[english]{babel}
\usepackage[latin1]{inputenc}
\usepackage[T1]{fontenc}
\usepackage{graphicx}
\usepackage{amsmath}

\usepackage{amssymb}
\usepackage{latexsym}
\usepackage{amsthm}
\usepackage[all]{xy}
\usepackage{stmaryrd}
\usepackage{subfigure}
\usepackage{MnSymbol}
\usepackage{comment} 

\def\be{\begin{equation}}
\def\ee{\end{equation}}
\def\beq{\begin{eqnarray*}}
\def\eeq{\end{eqnarray*}}
\def\Z{\mathbb{Z}}

\newcommand{\pic}[3]{\parbox[c]{#1cm}{\includegraphics[scale=#2]{#3}}}

\newtheorem{theo}{Theorem}[section]

\newtheorem{lem}[theo]{Lemma}
\newtheorem{prop}[theo]{Proposition}

\theoremstyle{definition}
\newtheorem{defn}[theo]{Definition}

\newtheorem{quest}[theo]{Question}
\newtheorem{rem}[theo]{Remark}

\author{Alessio Carrega}
\address{Dipartimento di Matematica, Largo Pontecorvo 5, 56127 Pisa, Italy}
\email{carrega at mail dot dm dot unipi dot it}

\title[$9$ generators of the skein space of the 3-torus]{$9$ generators of the skein space of the 3-torus}

\begin{document}

\begin{abstract}
We show that the skein vector space of the 3-torus is finitely generated. We show that it is generated by $9$ elements: the empty set, some simple closed curves representing the non null elements of the first homology group with coefficients in $\Z_2$, and a link consisting of two parallel copies of one of the previous non empty knots.
\end{abstract}

\maketitle

\setcounter{tocdepth}{1}

\section{Introduction}

An alternative approach to representation theory for \emph{quantum invariants} is provided by \emph{skein theory}. The word ``skein'' and the notion were introduced by Conway in 1970 for his model of the \emph{Alexander polynomial}. This idea became really useful after the work of Kauffman \cite{Kauffman} which redefined the \emph{Jones polynomial} in a very simple and combinatorial way passing through the \emph{Kauffman bracket}. These combinatorial techniques allow us to reproduce all quantum invariants arising from the representations of $U_q(\mathfrak{sl}_2)$ without any reference to representation theory. This also leads to many interesting and quite easy computations. This skein method was used by various authors \cite{Lickorish1, Lickorish2, Lickorish3, Lickorish4, BHMV, Kauffman-Lins} to re-interpret and extend some of the methods of representation theory.

The first notion in skein theory is the one of ``\emph{skein vector space}'' (or \emph{skein module}). These are vector spaces ($R$-modules) associated to oriented 3-manifolds, where the base field is equipped with a fixed invertible element $A$. These were introduced independently in 1988 by Turaev \cite{Turaev0} and in 1991 by Hoste and Przytycki \cite{HP0}. We can think of them as an attempt to get an algebraic topology for knots: they can be seen as homology spaces obtained using isotopy classes instead of homotopy or homology classes. In fact they are defined taking a vector space generated by sub-objects (\emph{framed links}) and then quotienting them by some relations. In this framework, the following questions arise naturally and are still open in
general:
\begin{quest}
$\ $
\begin{itemize}
\item{Are skein spaces (modules) computable?}
\item{How powerful are them to distinguish 3-manifolds and links?}
\item{Do the vector spaces (modules) reflect the topology/geometry of the 3-manifolds (\emph{e.g.} surfaces, geometric decomposition)?}
\item{Does this theory have a functorial aspect? Can it be extended to a functor from a category of cobordisms to the category of vector spaces (modules) and linear maps?}
\end{itemize}
\end{quest}

Skein spaces (modules) can also be seen as deformations of the ring of the $SL_2(\mathbb{C})$-\emph{character variety} of the 3-manifold \cite{Bullock2}. Moreover they are useful to generalize the Kauffman bracket, hence the Jones polynomial, to manifolds other than $S^3$. Thanks to result of Hoste-Przytycki \cite{HP2, Pr2} and (with different techniques) to Costantino \cite{Costantino2}, now we can define the Kauffman bracket also in the connected sum $\#_g(S^1\times S^2)$ of $g\geq 0$ copies of $S^1\times S^2$.

Until now there are only few 3-manifolds whose skein space (module) is known, see for instance \cite{Bullock0, HP1, HP3, HP2, Marche1, Mroczkowski1, Mroczkowski2, Mroczkowski-Dabkowski, Pr1, Pr2, Pr:survey, ThangLe}. Another natural question is: 
\begin{quest}
Is the skein vector space of a closed oriented 3-manifold always finite generated?
\end{quest}

In this paper we take as base field the set $\mathbb{Q}(A)$ of all rational functions with rational coefficients and abstract variable $A$, and then we note that every result in this work holds also for the field $\mathbb{C}$ of complex numbers with $A\in \mathbb{C}$ a non null number such that $A^n \neq 1$ for every $n>0$.

\begin{theo}
The skein space $K(T^3)$ of the 3-torus $T^3=S^1\times S^1\times S^1$ is finitely generated.
\end{theo}

A set of $9$ generators is given by the empty set $\varnothing$, some simple closed curves representing the non null elements of the first homology group $H_1(T^3;\Z_2)\cong (\Z_2)^3$ with coefficients in $\Z_2$, and a skein element $\alpha$ that is equal to the link consisting of two parallel copies of any previous non empty knots. 

Our main tool is the algebraic work of Frohman and Gelca \cite{Frohman-Gelca}. The skein space (module) of a (thickened) surface has a natural structure of algebra obtained by overlap of framed links. In their work Frohman and Gelca gave a nice formula that describes the product in the skein space (algebra) $K(T^2)$ of the 2-torus $T^2=S^1\times S^1$. A standard embedding of $T^2$ in $T^3$ makes this product commutative, hence we can get further relations from the formula of Frohman-Gelca.

\subsection*{Acknowledgments}

The author is warmly grateful to Bruno Martelli for his constant support and encouragement.

\section{The result}

\subsection{Definition of skein module}

Let $M$ be an oriented 3-manifold, $R$ a commutative ring with unit and $A\in R$ an invertible element of $R$. Let $V$ be the abstract free $R$-module generated by all framed links in $M$ (considered up to isotopies) including the empty set $\varnothing$. 

\begin{defn}
The $(R,A)$-\emph{Kauffman bracket skein module} of $M$, or the $R$-\emph{skein module}, or simply the \emph{KBSM}, is sometimes indicated with $KM(M;R,A)$, and is the quotient of $V$ by all the possible \emph{skein relations}:
$$
\begin{array}{rcl}
 \pic{1.2}{0.3}{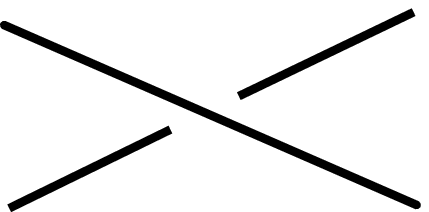}  & = & A \pic{1.2}{0.3}{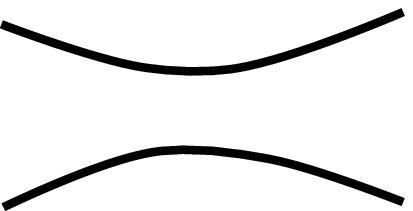}  + A^{-1}  \pic{1.2}{0.3}{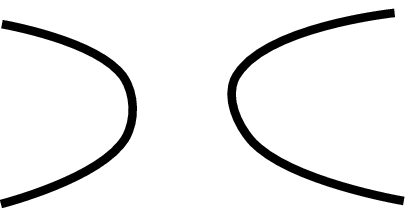}  \\
 L \sqcup \pic{0.8}{0.3}{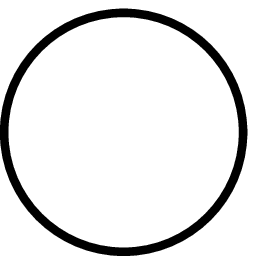}  & = & (-A^2 - A^{-2})  D  \\
\pic{0.8}{0.3}{banp.eps}  & = & (-A^2-A^{-2}) \varnothing 
\end{array} .
$$
These are local relations where the framed links in an equation differ just in the pictured 3-ball that is equipped with a positive trivialization. An element of $KM(M;R,A)$ is called a \emph{skein} or a \emph{skein element}. If $M$ is the oriented $I$-bundle over a surface $S$ (this is $M=S\times [-1,1]$ if $S$ is oriented) we simply write $KM(S;R,A)$ and call it the \emph{skein module} of $S$.

Let $\mathbb{Q}(A)$ be field of all rational function with rational coefficients and abstract variable $A$. We set
$$
K(M) := KM(M; \mathbb{Q}(A) ,A) 
$$
and we call it the \emph{skein vector space}, or simply the \emph{skein space}, of $M$.
\end{defn}

\begin{rem}
It is easy to verify that if we modify the framing of a component of a framed link, the skein changes by the multiplication of an integer power of $-A^3$: 
\beq
\pic{1.2}{0.3}{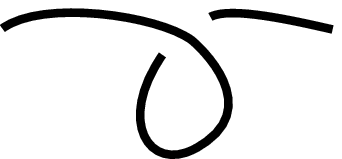} & = & -A^3 \pic{1.2}{0.3}{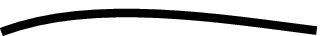} \\
\pic{1.2}{0.3}{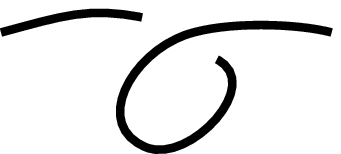} & = & -A^{-3} \pic{1.2}{0.3}{riga.eps} .
\eeq
\end{rem}

\subsection{The skein algebra of the 2-torus}

\begin{defn}
Let $S$ be a surface the skein module $KM(S;R,A)$ has a natural structure of $R$-algebra that is given by the linear extension of the multiplication defined on framed links. Given two framed links $L_1,L_2 \subset S\times [-1,1]$, the product $L_1\cdot L_2 \subset S\times [-1,1]$ is obtained by putting $L_1$ above $L_2$, $L_1\cdot L_2 \cap S\times [0,1] = L_1$ and $L_1\cdot L_2 \cap S\times [-1,0] = L_2$.
\end{defn}

Look at the 2-torus $T^2$ as the quotient of $\mathbb{R}^2$ modulo the standard lattice of translations generated by $(1,0)$ and $(0,1)$, hence for any non null pair $(p,q)$ of integers we have the notion of $(p,q)$-\emph{curve}: the simple closed curve in the 2-torus that is the quotient of the line passing trough $(0,0)$ and $(p,q)$.

\begin{defn}
Let $p$ and $q$ be two co-prime integers, hence $(p,q)\neq (0,0)$. We denote by $(p,q)_T$ the $(p,q)$-curve in the 2-torus $T^2$ equipped with the black-board framing. Given a framed knot $\gamma$ in an oriented 3-manifold $M$ and an integer $n\geq 0$, we denote by $T_n(\gamma)$ the skein element defined by induction as follows:
\beq
T_0(\gamma) & := & 2\cdot \varnothing  \\
T_1(\gamma) & := & \gamma \\
T_{n+1}(\gamma) & := & \gamma \cdot T_n(\gamma) - T_{n-1}(\gamma)
\eeq
where $\gamma \cdot T_n(\gamma)$ is the skein element obtained adding a copy of $\gamma$ to all the framed links that compose the skein $T_n(\gamma)$. For $p,q\in \Z$ such that $(p,q)\neq (0,0)$, we denote by $(p,q)_T$ the skein element 
$$
(p,q)_T := T_{{\rm MCD}(p,q)} \left( \left( \frac{p}{ {\rm MCD}(p,q) } , \frac{q}{ {\rm MCD}(p,q) } \right)_T \right) ,
$$
where ${\rm MCD}(p,q)$ is the maximum common divisor of $p$ and $q$. Finally we set 
$$
(0,0)_T := 2 \cdot \varnothing .
$$
\end{defn}

It is easy to show that the set of all the skein elements $(p,q)_T$ with $p,q\in\Z$ generates $KM(T^2;R,A)$ as $R$-module.

This is not the standard way to color framed links in a skein module. The colorings $JW_n(\gamma)$, $n\geq 0$, with the Jones-Wenzl projectors are defined in the same way of $T_n(\gamma)$ but at the $0$-level we have $JW_0(\gamma)=\varnothing$.

\begin{theo}[Frohman-Gelca]\label{theorem:F-G}
For any $p,q,r,s\in \Z$ the following holds in the skein module $KM(T^2;R,A)$ of the 2-torus $T^2$:
$$
(p,q)_T \cdot (r,s)_T = A^{\left| \begin{matrix} p & q \\ r & s \end{matrix} \right|} (p+r,q+s)_T + A^{-\left| \begin{matrix} p & q \\ r & s \end{matrix} \right|} (p-r,q-s)_T ,
$$
where $\left| \begin{matrix} p & q \\ r & s \end{matrix} \right|$ is the determinant $ps - qr$.
\begin{proof}
See \cite{Frohman-Gelca}.
\end{proof}
\end{theo}

\subsection{The abelianization}

\begin{defn}
Let $B$ be a $R$-algebra for a commutative ring with unity $R$. We denote by $C(B)$ the $R$-module defined as the following quotient: 
$$
C(B) := \frac{B}{ [ B , B ] }  
$$
where $[B , B ]$ is the sub-module of $B$ generated by all the elements of the form $ab- ba$ for $a,b\in B$. We call $C(B)$ the \emph{abelianization} of $B$.
\end{defn}

\begin{rem}
Usually in non-commutative algebra the \emph{abelianization} is the $R$-algebra defined as the quotient of $B$ modulo the sub-algebra (sub-module and ideal) generated by all the elements of the form $ab-ba$. In our definition the denominator is just a sub-module and we only get a $R$-module. We use the word ``abelianization'' anyway. 
\end{rem}

Now we work with $C(K(T^2))$ and we still use $(p,q)_T$ and $(p,q)_T \cdot (r,s)_T$ to denote the class of $(p,q)_T\in K(T^2)$ and $(p,q)_T\cdot (r,s)_T\in K(T^2)$ in $C(K(T^2))$.

\begin{lem}\label{lem:ab_alg_2-tor}
Let $(p,q)$ be a pair of integers different from $(0,0)$. Then in the abelianization $C(K(T^2))$ of the skein algebra $K(T^2)$ of the 2-torus $T^2$ we have
$$
(p,q)_T = \begin{cases}
(1,0)_T & \text{if } p\in 2\Z+1 , \ q \in 2\Z \\
(0,1)_T & \text{if } p\in 2\Z , \ q \in 2\Z+1 \\
(1,1)_T & \text{if } p,q\in 2\Z+1 \\
(2,0)_T & \text{if } p,q \in 2\Z
\end{cases} .
$$
Hence $C(K(T^2))$ is generated as a $\mathbb{Q}(A)$-vector space by the empty set $\varnothing$, the framed knots $(1,0)_T$, $(0,1)_T$, $(1,1)_T$, and a framed link consisting of two parallel copies of $(1,0)_T$.
\begin{proof}
By Theorem~\ref{theorem:F-G} for every $p,q\in \Z$ we have
\beq
A^{-q} (p+2,q)_T + A^q (p,q)_T & = & (p+1,q)_T \cdot (1,0)_T \\
 & = & (1,0)_T \cdot (p+1,q)_T \\
 & = & A^q (p+2,q)_T + A^{-q} (-p,-q)_T .
\eeq
Since $(p,q)_T=(-p,-q)_T$ we have $ (A^q - A^{-q}) (p,q)_T = (A^q -A^{-q}) (p+2,q)_T $. Hence if $q\neq 0$ we get $(p,q)_T = (p+2,q)_T $ (here we use the fact that the base ring is a field and $A^n \neq 1$ for every $n>0$). Thus
$$
(p,q)_T = \begin{cases}
(0,q)_T & \text{if } p\in 2\Z, \ \ q\neq 0 \\ 
(1,q)_T & \text{if } p\in 2\Z +1, \ \ q\neq 0
\end{cases} .
$$
Analogously by using $(0,1)_T$ instead of $(1,0)_T$ for $p\neq 0$ we get
$$
(p,q)_T = \begin{cases}
(p,0)_T & \text{if } q\in 2\Z, \ \ q\neq 0 \\ 
(p,1)_T & \text{if } q\in 2\Z +1, \ \ q\neq 0
\end{cases} .
$$
Therefore if $p,q\in 2\Z+1$, $(p,q)_T= (1,1)_T$. If $p\neq 0$ we get
$$
(p,0)_T=(p,2)_T = \begin{cases}
(0,2)_T & \text{if } p\in 2\Z \\
(1,2)_T & \text{if } p\in 2\Z +1
\end{cases} 
= \begin{cases}
(0,2)_T & \text{if } p\in 2\Z \\
(1,0)_T & \text{if } p\in 2\Z +1
\end{cases} .
$$
In the same way for $q\neq 0$ we get
$$
(0,q)_T = (2,q)_T = \begin{cases}
(2,0)_T & \text{if } p\in 2\Z \\
(2,1)_T & \text{if } p\in 2\Z +1
\end{cases} 
= \begin{cases}
(2,0)_T & \text{if } p\in 2\Z \\
(0,1)_T & \text{if } p\in 2\Z +1
\end{cases} .
$$
In particular we have
$$
(2,0)_T = (2,2)_T = (2,-2)_T = (0,2)_T = (p,q)_T  \text{ for } (p,q)\neq (0,0), \ p,q \in 2\Z .
$$
\end{proof}
\end{lem}

\subsection{The $(p,q,r)$-type curves}

As for the 2-torus $T^2$, we look at the 3-torus $T^3$ as the quotient of $\mathbb{R}^3$ modulo the standard lattice of translations generated by $(1,0,0)$, $(0,1,0)$ and $(0,0,1)$.

\begin{defn}
Let $(p,q,r)$ be a triple of co-prime integers, that means ${\rm MCD}(p,q,r)=1$, where ${\rm MCD}(p,q,r)$ is the maximum common divisor of $p$, $q$ and $r$, in particular we have $(p,q,r)\neq (0,0,0)$. The $(p,q,r)$-\emph{curve} is the simple closed curve in the 3-torus that is the quotient (under the standard lattice) of the line passing through $(0,0,0)$ and $(p,q,r)$. We denote by $[p,q,r]$ the $(p,q,r)$-curve equipped with the framing that is the collar of the curve in the quotient of any plane containing $(0,0,0)$ and $(p,q,r)$. The framing does not depend on the choice of the plane.
\end{defn}

\begin{defn}
An embedding $e:T^2\rightarrow T^3$ of the 2-torus in the 3-torus is \emph{standard} if it is the quotient (under the standard lattice) of a plane in $\mathbb{R}^3$ that is the image of the plane generated by $(1,0,0)$ and $(0,1,0)$ under a linear map defined by a matrix of $SL_3(\Z)$ (a $3\times 3$ matrix with integer entries and determinant $1$).
\end{defn}

\begin{rem}\label{rem:std_emb}
There are infinitely many standard embeddings even up to isotopies. A standard embedding of $T^2$ in $T^3$ is the quotient under the standard lattice of the plane generated by two columns of a matrix of $SL_3(\Z)$.
\end{rem}

\begin{lem}\label{lem:pqr-curve}
Let $(p,q,r)$ be a triple of co-prime integers. Then the skein element $[p,q,r]\in K(T^3)$ is equal to $[x,y,z]$, where $x,y,z\in \{0,1\}$ and have respectively the same parity of $p$, $q$ and $r$.
\begin{proof}
Every embedding $e:T^2 \rightarrow T^3$ of the 2-torus $T^2$ in $T^3$ defines a linear map between the skein spaces
$$
e_* : K(T^2) \longrightarrow K(T^3) .
$$
The map $e_*$ factorizes with the quotient map $K(T^2) \rightarrow C(K(T^2))$. In fact we can slide the framed links in $e(T^2\times [-1,1])$ from above to below getting $e_*(L_1\cdot L_2) = e_*(L_2\cdot L_1)$ for every two framed links, $L_1$ and $L_2$, in $T^2\times [-1,1]$. As said in Remark~\ref{rem:std_emb}, a standard embedding $e:T^2 \rightarrow T^3$ corresponds to the plane generated by two columns $(p_1,q_1,r_1) , (p_2,q_2,r_2) \in \Z^3$ of a matrix in $SL_3(\Z)$. In this correspondence $e_*(( a,b)_T)= [ap_1 + bp_2, aq_1 +b q_2, ar_1 + br_2]$ for every co-prime $a,b \in \Z$. Therefore by Lemma~\ref{lem:ab_alg_2-tor} we get 
\beq
[a'p_1 + b'p_2, a'q_1 +b' q_2, a'r_1 + b'r_2] & = & e_*(( a',b')_T) \\
& = & e_*(( a,b)_T) \\
& = & [ap_1 + bp_2, aq_1 +b q_2, ar_1 + br_2]
\eeq
for every two pairs $(a,b),(a',b')\in \Z^2$ of co-prime integers such that $a+a', b+b'\in 2\Z$.

Let $(p,q,r)$ be a tripe of co-prime integers. By permuting $p,q,r$ we get either $(p,q,r)=(1,0,0)$ or $p,q\neq 0$. Consider the case $p,q\neq 0$. Let $d$ be the maximum common divisor of $p$ and $q$, and let $\lambda,\mu \in \Z$ such that $\lambda p + \mu q = d$. The following matrix belongs in $SL_3(\Z)$:
$$
M_1:= \left( \begin{matrix}
\frac p d & -\mu & 0 \\
\frac q d & \lambda & 0 \\
0 & 0 & 1
\end{matrix} \right) .
$$
Let $v^{(1)}_1$ and $v^{(1)}_3$ be the first and the third columns of $M_1$. We have $(p,q,r)= dv^{(1)}_1 + r v^{(1)}_3$. Hence
$$
[p,q,r] = \begin{cases}
[\frac p d , \frac q d, 0]  & \text{if } d\in 2\Z+1 ,\ r\in 2\Z \\
[0 , 0 , 1]  & \text{if } d \in 2\Z, \  r\in 2\Z+1 \\
[\frac p d , \frac q d, 1] & \text{if } d,r\in 2\Z+1 
\end{cases} .
$$
The integers $p,q,r$ can not be all even because they are co-prime, hence $d$ and $r$ can not be both even. Therefore we just need to study the cases where $r\in \{0,1\}$. 

If $r=0$ we consider the trivial embedding of $T^2$ in $T^3$. The corresponding matrix of $SL_3(\Z)$ is the identity. We have $(p/d,q/d,0)=p/d(1,0,0) + q/d(0,1,0)$, hence
$$
[p,q,0] = [\frac p d , \frac q d ,0] = \begin{cases}
[1,0,0] & \text{ if } \frac p d \in 2\Z+1, \ \frac q d \in 2\Z \\
[0,1,0] & \text{ if } \frac p d \in 2\Z, \ \frac q d \in 2\Z+1 \\
[1,1,0] & \text{ if } \frac p d , \frac q d \in 2\Z+1 
\end{cases} .
$$

If $r=1$ we take the following matrix of $SL_3(\Z)$:
$$
M_2:= \left( \begin{matrix}
0 & 0 & 1 \\
q & -1 & 0 \\
1 & 0 & 0
\end{matrix} \right) .
$$
Let $v^{(2)}_1$ and $v^{(2)}_3$ be the first and the third columns of $M_2$. We have $(p,q,1)= pv^{(2)}_3 + v^{(2)}_1$, hence
$$
[p,q,1] = \begin{cases}
[1,q,1] & \text{if } p\in 2\Z+1 \\
[0,q,1] & \text{if } p\in 2\Z 
\end{cases} .
$$

By permuting $p,q,r$ we reduce the case $(p,q,r)=(0,q,1)$ to the case $p,q\neq 0$, $r=0$ that we studied before.
 
It remains to consider the case $p=r=1$. We consider the following matrix of $SL_3(\Z)$:
$$
M_3:= \left( \begin{matrix}
1 & 0 & 0 \\
0 & 1 & 0 \\
1 & 0 & 1
\end{matrix} \right) .
$$
Let $v^{(3)}_1$ and $v^{(3)}_2$ be the first and the second columns of $M_3$. We have $(1,q,1)= v^{(3)}_1 + qv^{(3)}_2$. Hence
$$
[1,q,1] = \begin{cases}
[1,0,1] & \text{if } q\in 2\Z \\
[1,1,1] & \text{if } q\in 2\Z+1 
\end{cases} .
$$
\end{proof}
\end{lem}

\begin{lem}\label{lem:intersect_2-tor}
The intersection of any two different standard embedded 2-tori in $T^3$ contains a $(p,q,r)$-type curve.
\begin{proof}
Let $T_1$ and $T_2$ be two standard embedded tori in the 3-torus, and let $\pi_1$ and $\pi_2$ be two planes in $\mathbb{R}^3$ whose projection under the standard lattice is respectively $T_1$ and $T_2$. The intersection $T_1\cap T_2$ contains the projection of $\pi_1\cap \pi_2$. We just need to prove that in $\pi_1\cap \pi_2$ there is a point $(p,q,r)\neq (0,0,0)$ with integer coordinates $p,q,r\in \Z$. Every plane defining a standard embedded torus is generated by two vectors with integer coordinates, and hence it is described by an equation $ax+by+cz=0$ with integer coefficients $a,b,c\in \Z$. Applying a linear map described by a matrix of $SL_3(\Z)$ we can suppose that $\pi_1$ is the trivial plane $\{ z=0 \}$. Let $a,b,c\in \Z$ such that $\pi_2= \{ax+by+cz=0 \}$. The vector $(-b,a,0)$ is non null and lies on $\pi_1\cap \pi_2$.
\end{proof}
\end{lem}

\subsection{Diagrams}

\begin{figure}[htbp]
\begin{center}
\end{center}
\includegraphics[width = 9cm]{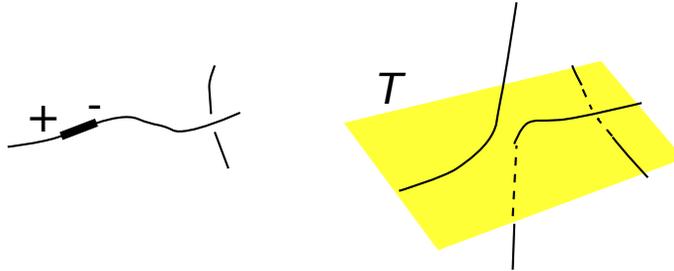}
\caption{Diagrams of framed links in $T^3$. The yellow plane is a part of the standard embedded torus $T\subset T^3$ where the links project. If we look at the framed links in $T^3$ as framed tangles in $T\times [-1,1]$, the two strands that get out vertically from the yellow plane end in the boundary points $(x,1)$ and $(x,-1)$ for some $x\in T$.}
\label{figure:diag_3-tor}
\end{figure}

Framed links in $T^3$ can be represented by diagrams in the 2-torus $T^2$. These diagrams are like the usual link diagrams but with further oriented signs on the edges (see Fig.~\ref{figure:diag_3-tor}-(left)). Fix a standard embedded 2-torus $T$ in $T^3$. After a cut along a parallel copy $T'$ of $T$, the 3-torus becomes diffeomorphic to $T\times [-1,1]$ and framed links in $T^3$ correspond to framed tangles of $T\times [-1,1]$. These diagrams are generic projections on $T$ of the framed tangles in $T\times [-1,1]$ via the natural projection $(x,t)\mapsto x$. The further signs on the diagrams represent the intersection of the framed links with the boundary $T\times \{-1,1\}$, in other words they represent the passages of the links along the $(p,q,r)$-type curve that in the Euclidean metric is orthogonal to $T$ (see Fig.~\ref{figure:diag_3-tor}-(right)). If $T$ is the trivial torus $S^1\times S^1 \times \{x\}$, the further signs represent the passages through the third $S^1$-factor. We use the proper notion of blackboard framing.

\subsection{Generators for the 3-torus}

The following is the main theorem proved in this paper. We use all the previous lemmas to get a set of $9$ generators of $K(T^3)$.

\begin{theo}\label{theorem:sk_3-tor}
The skein space $K(T^3)$ of the 3-torus $T^3$ is generated by the empty set $\varnothing$, $[1,0,0]$, $[0,1,0]$, $[0,0,1]$, $[1,1,0]$, $[1,0,1]$, $[0,1,1]$, $[1,1,1]$ and a skein $\alpha$ that is equal to the framed link consisting of two parallel copies of any $(p,q,r)$-type curve.
\begin{proof}
Let $T$ be the trivial embedded 2-torus: the one containing the $(p,q,r)$-type curves with $r=0$. Use $T$ to project the framed links and make diagrams. By using the first skein relation on these diagrams we can see that $K(T^3)$ is generated by the framed links described by diagrams on $T$ without crossings. These diagrams are union of simple closed curves on $T$ equipped with some signs as the one with $+$ and $-$ in Fig.~\ref{figure:diag_3-tor}. These simple closed curves are either parallel to a $(p,q)$-curve or homotopically trivial. The framed links described by these diagrams lie in the standard embedded tori that are the projection (under the standard lattice) of the planes generated by $(0,0,1)$ and $(p,q,0)$ for some $p$ and $q$. Therefore $K(T^3)$ is generated by the images of $K(T^2)$ under the linear maps induced by the standard embeddings of $T^2$ in $T^3$.

As said in the proof of Lemma~\ref{lem:pqr-curve}, the linear map $e_*$ induced by any standard embedding $e:T^2\rightarrow T^3$ factorizes with the quotient map $K(T^2) \rightarrow C( K(T^2) )$. Lemma~\ref{lem:ab_alg_2-tor} applied on the standard embedding $e$ shows that the image $e_*(K(T^2))$ is generated by $\varnothing$, three $(p,q,r)$-type curves lying on $e(T^2)$, and the skein $\alpha_e$ that is equal to the framed link consisting of two parallel copies of any $(p,q,r)$-type curve lying on $e(T^2)$.

Let $e_1,e_2:T^2 \rightarrow T^3$ be two standard embeddings. By Lemma~\ref{lem:intersect_2-tor} $e_1(T^2)\cap e_2(T^2)$ contains a $(p,q,r)$-type curve $\gamma$, hence $\alpha_{e_1}$ and $\alpha_{e_2}$ coincide with the framed link that is two parallel copies of $\gamma$. Therefore the skein element $\alpha_e$ does not depend on the embedding $e$.

We conclude by using Lemma~\ref{lem:pqr-curve} that says that the skein of any $(p,q,r)$-type curve is equal to the one of a standard representative of a non null element of the first homology group $H_1(T^3;\Z_2)$ with coefficient in $\Z_2$, namely a $(p,q,r)$-type curve with $p,q,r\in \{0,1\}$.
\end{proof}
\end{theo}

\begin{rem}
Theorem~\ref{theorem:sk_3-tor}, Lemma~\ref{lem:ab_alg_2-tor} and Lemma~\ref{lem:pqr-curve} work for every base pair $(R,A)$ such that $A^n-1$ is an invertible element of $R$ for any $n>0$. In particular they work for $(\mathbb{C}, A)$, where $A^n\neq 1$ for any $n>0$.
\end{rem}

\subsection{Linear independence}

Here we talk about the linear independence of generators of $K(T^2)$ we have shown. The following proposition shows a decomposition in direct sum of $K(T^3)$.

\begin{prop}\label{prop:3-tor}
The skein space $K(T^3)$ is the direct sum of $8$ sub-spaces
$$
K(T^3) = V_0 \oplus V_1 \oplus \ldots \oplus V_7 
$$
such that:
\begin{enumerate}
\item{$V_0$ is generated by the empty set $\varnothing$ and the skein $\alpha$ (see Theorem~\ref{theorem:sk_3-tor});}
\item{every $(p,q,r)$-type curve generates a $V_j$ with $j>0$ and
every $V_j$ with $j>0$ is generated by one such curve.}
\end{enumerate}
\begin{proof}
The skein relations relates framed links holding in the same $\Z_2$-homology class. Hence for every oriented 3-manifold $M$ we have a decomposition in direct sum 
$$
KM(M;R,A) = \bigoplus_{h \in H_1(M;\Z_2)} V_h ,
$$
where $V_h$ is generated by the framed links whose $\Z_2$-homology class is $h$. The statement follow by this observation and the fact that if $[p,q,r]$ and $[p',q',r']$ represent the same $\Z_2$-homology class, then $[p,q,r]=[p',q',r']\in K(T^3)$.
\end{proof}
\end{prop}

\begin{rem}
Given a triple of integers $(x,y,z)\neq (0,0,0)$ such that $x,y,z\in \{0,1\}$, we can easily find an orientation preserving diffeomorphism of the 3-torus $T^3$ sending $[x,y,z]$ to $[1,0,0]$. Hence if the skein of one such curve $[x,y,z]$ is null then also all the others skein elements of such curves are null. Therefore by Proposition~\ref{prop:3-tor} the possibly dimensions of the skein space $K(T^3)$ are $0$, $1$, $2$, $7$, $8$ and $9$.
\end{rem}

\begin{quest}
Is $9$ the dimension of the skein vector space $K(T^3)$ of the 3-torus?
\end{quest}

\end{document}